\documentclass[12pt]{article}
\usepackage{amsmath,amsthm,amssymb}
\usepackage{mathrsfs}
\usepackage[english]{babel}

\newtheorem{theorem}{Theorem}[section]

\newtheorem{lemma}{Lemma}[section]

\theoremstyle{definition}

\theoremstyle{remark}
\newtheorem{remark}{Remark}[section]

\DeclareMathOperator{\dom}{\operatorname{Dom}}

\newcommand\R{{\mathbb R}}
\newcommand\X{{\R^d}}
\newcommand\N{{\mathbb N}}

\newcommand\M{{\mathcal M}}
\newcommand\F{{\mathcal F}}
\newcommand\EE{{\mathcal E}}

\newcommand\G{{\mathcal G}}

\newcommand\B{{\mathcal B}}
\newcommand\La{\Lambda}

\newcommand\Ga{\Gamma}
\newcommand\ga{\gamma}
\newcommand\GL{{\Ga_\La}}

\newcommand\eps{\varepsilon}

\newcommand\lmod{\left\vert}
\newcommand\rmod{\right\vert}
\newcommand{\sk}[1]{{\left( #1 \right)}}
\newcommand{\skf}[1]{{\left\{ #1 \right\}}}
\newcommand{\skp}[1]{{\left[ #1 \right]}}
\newcommand{\skg}[1]{{\left\langle #1 \right\rangle}}
\newcommand{\md}[1]{{\lmod #1 \rmod}}

\newcommand\goto{\rightarrow}
\allowdisplaybreaks[3]

\numberwithin{equation}{section}

\newcommand{\iX}{\int_{\X}}
\newcommand{\iXn}[1]{\int_{{\sk{\X}}^{#1}}}
\newcommand{\iXg}{\int_{\X}\ga\sk{dx}}
\newcommand{\iXx}{\iX z\, dx}
\newcommand{\iXy}{\iX dy}
\newcommand{\iG}{\int_{\Ga}}
\newcommand{\iGG}{\int_{\Ga}\mu\sk{d\ga}}
\newcommand{\aeps}[1]{ a_\varepsilon(#1)}
\newcommand{\f}[1]{\sk{e^{\varphi\sk{#1}}-1}}
\newcommand{\g}[1]{e^{\varphi\sk{#1}}}
\newcommand{\vect}[2]{\sk{#1_1,\ldots,#1_{#2}}}
\newcommand{\dvect}[2]{\,d #1_1 \dotsm d#1_{#2}}

\newcommand{\rom}[1]{{\rm #1}}

\begin{document}

\begin{center}{\Large \bf
  Equilibrium Glauber dynamics\\ of continuous particle systems\\ as a scaling limit of  Kawasaki dynamics
}\end{center}

{\large Dmitri  L. Finkelshtein}\\
Institute of Mathematics, National Academy of Sciences of Ukraine, 3 Tereshchenkivska Str., 
Kiev 01601, Ukaine.\\
e.mail: \texttt{fdl@imath.kiev.ua}\vspace{2mm}

{\large Yuri G. Kondratiev}\\
 Fakult\"at f\"ur Mathematik, Universit\"at
Bielefeld, Postfach 10 01 31, D-33501 Bielefeld, Germany; BiBoS, Univ.\ Bielefeld, Germany; Kiev-Mohyla Academy, Kiev, Ukraine.\\
 e-mail:
\texttt{kondrat@mathematik.uni-bielefeld.de}\vspace{2mm}

{\large Eugene W. Lytvynov}\\
Department of Mathematics, University of Wales Swansea, Singleton
Park, Swansea SA2 8PP, U.K.\\ e-mail:
\texttt{e.lytvynov@swansea.ac.uk}

{\small
\begin{center}
{\bf Abstract}
\end{center}
\noindent A Kawasaki dynamics in continuum is a dynamics of an infinite system 
of interacting particles in $\R^d$ which randomly hop over the space. 
In this paper, we deal with an equilibrium Kawasaki dynamics which has a Gibbs measure $\mu$ as  
invariant measure. We study a scaling  limit  
of such a dynamics, derived through a scaling of the jump rate. Informally, we expect that, in the limit, only jumps of ``infinite length'' will survive, i.e., we expect to arrive at a Glauber dynamics in continuum (a birth-and-death process in $\R^d$). We prove that, in the low activity-high temperature 
regime, the generators of the Kawasaki dynamics converge to the generator of a 
Glauber dynamics. The convergence is on the set of exponential functions, in the $L^2(\mu)$-norm. Furthermore, additionally assuming that the potential of pair interaction
is positive, we prove the weak convergence of the finite-dimensional distributions of the processes.  
 }
\vspace{3mm}

\noindent 
{\it MSC:} 60K35, 60J75, 60J80, 82C21, 82C22 \vspace{1.5mm}

\noindent{\it Keywords:} Continuous
system; Gibbs measure; Glauber dynamics; Kawasaki dynamics; Scaling limit\vspace{1.5mm}


\section{Introduction} A Kawasaki dynamics in continuum is a dynamics of an infinite system 
of interacting particles in $\R^d$ which randomly hop over the space. In this paper, we deal with an equilibrium 
Kawasaki dynamics which has a Gibbs measure $\mu$ as invariant measure. About  $\mu$ we assume that it corresponds to an activity parameter $z>0$ and a potential of pair interaction $\phi$. The generator of the Kawasaki dynamics 
is given, on an appropriate set of cylinder functions, by
\begin{align} (HF)(\gamma)&=-\sum_{x\in\gamma}\int_{\R^d}dy\, a(x-y)\exp\left[-\sum_{u\in\gamma\setminus x}\phi(u-y)\right]
\notag\\ &\qquad\times (F(\gamma\setminus x\cup y)-F(\gamma)),\quad \gamma\in\Gamma.\label{fff}\end{align} 
Here, $\Gamma$ denotes the configuration space over $\R^d$, i.e., the space of all locally finite subsets of $\R^d$,
and, for simplicity of notations, we just write $x$ instead of $\{x\}$. About the function $a(\cdot)$ in \eqref{fff} we assume that it is non-negative, 
integrable and symmetric with respect to the origin. 
 The factor $a(x-y)\exp\left[-\sum_{u\in\gamma\setminus x}\phi(u-y)\right]$ in \eqref{fff}
describes the rate with which, given a configuration $\gamma\in\Gamma$, a particle $x\in\gamma$ jumps to $y$.

Under very mild assumptions on the Gibbs measure $\mu$, it was proved in \cite{KLR} that there indeed exists a Markov process on $\Gamma$ with {\it c\'adl\'ag} paths whose generator is given by
\eqref{fff}. We assume that the initial distribution of this dynamics is $\mu$, and perform the following scaling
of this dynamics. For each $\eps>0$, we consider the equilibrium Kawasaki dynamics whose generator is given by formula \eqref{fff} in which $a(\cdot)$ is replace by the function \begin{equation}\label{rtsra} a_\eps(\cdot):=\eps^da(\eps\cdot).\end{equation}  
We denote this generator by $H_\eps$, and study the limit of the corresponding dynamics as $\varepsilon\to0$.
We would like to stress that we scale the jump rate of each particle, whereas the interaction between particles remains the same. It means that such a scaling is not at all of the mean-field type, as it might seem to be.

Informally, we expect that, in the limit, only jumps of infinite length will survive, i.e., jumps from a point to `infinity' and  from `infinity'
to a point. Thus, we expect to arrive at a Glauber dynamics in continuum, i.e., a birth-and-death process
in $\R^d$, cf.\ \cite{KL,KLR}. In fact, heuristic calculations show that the limiting Glauber dynamics has the generator 
\begin{align}
(H_0F)(\gamma)&=-\alpha\sum_{x\in\gamma}(F(\gamma\setminus x)-F(\gamma))\notag\\ &\quad -\alpha\int_{\R^d}z\,dx\,
\exp\left[-\sum_{u\in\gamma}\phi(u-x)\right] (F(\gamma\cup x)-F(\gamma)),\label{uigiy}
\end{align}
where 
\begin{equation}\label{iyuyd}
\alpha=z^{-1}k_\mu^{\sk{1}}\int_\X a\sk{x}\,dx,\end{equation}
$k_\mu^{(1)}$ being the first correlation function of the measure $\mu$. Thus, $\alpha$ describes the 
rate with  which a particle $x\in\gamma$ dies, whereas $\alpha z\exp\left[-\sum_{u\in\gamma}\phi(u-x)\right]$
describes the rate with which, given a configuration $\gamma$, a new particle is born at $x\in\R^d\setminus \gamma$.    
The existence of a Markov process on $\Gamma$ with {\it c\'adl\'ag} paths, whose generator is given by \eqref{uigiy}, was proved in \cite{KL} (see also \cite{KLR}). 

The main results of this paper are as follows:
\begin{itemize}

\item For any stable potential $\phi$ in the low activity-high temperature regime, the generators $H_\eps$
converge to the generator $H_0$. The convergence is on the set of exponential functions, in the $L^2(\Gamma, \mu)$-norm. 

\item For any positive potential $\phi$  in the low activity-high temperature regime, the finite-dimensional distributions
of the  Kawasaki dynamics with generator $H_\eps$ and initial distribution $\mu$ weakly converge to the finite-dimensional distributions of the Glauber dynamics with generator $H_0$ and initial distribution $\mu$. 

\end{itemize}

To prove the first main result, we essentially use the Ruelle bound on the correlation functions of the measure $\mu$, as well as the integrability of the Ursell (cluster) functions of $\mu$, see \cite{wewewe}. To derive from here the convergence of the 
finite-dimensional distributions  of the dynamics, we additionally need that the set of finite sums of exponential functions forms a core for the generator of the limiting dynamics, $H_0$. For this, we use a result from \cite{KL} on a core for $H_0$, which holds under the assumptions  
of positivity of the potential $\phi$.

We note that the generator of the Kawasaki dynamics is independent 
of the activity parameter $z>0$. Hence, at least heuristically, 
the Kawasaki dynamics has a continuum of symmetrizing Gibbs measures, indexed by the activity $z>0$. 
On the other hand, the limiting Glauber dynamics has only one of these measures as the symmetrizing one.  Thus, the result of the scaling  essentially depends on the initial distribution of the dynamics.

The paper is organized as follows. In Section \ref{iu}, we recall some known facts about Gibbs measures on the configuration space $\Gamma$. In Section \ref{uitgyuiftu}, we recall a rigorous construction of the equilibrium Kawasaki and Glauber dynamics. Our two main results are proved in Sections  \ref{section4} and \ref{section5}, respectively. Finally, in Section \ref{section6}, we make remarks on the results obtained, and  discuss some related open problems.

\section{Gibbs measures in the low activity-high temperature regime}
\label{iu}

The configuration space over $\X$,  $d\in\N$,
is defined by 
\[
\Ga:=\{\ga\subset\X:\, \md{\ga_\La}<\infty \text{
for each compact } \La\subset\X \},
\]
where $\md{\cdot}$ denotes the cardinality of a set and
$\ga_\La:=\ga\cap\La$. One can identify any $\ga\in\Ga$
with the positive Radon measure $\sum_{x\in\ga}\eps_x \in
\M\sk{\X}$, where $\eps_x$ is the Dirac measure with mass
at $x$, $\sum_{x\in\varnothing}\eps_x:=$zero measure, and
$\M\sk{\X}$ stands for the set of all positive Radon measures
on the Borel $\sigma$-algebra $\B\sk{\X}$. The space $\Ga$
can be endowed with the relative topology as a subset of
the space $\M\sk{\X}$ with the vague topology, i.e., the
weakest topology on $\Ga$ with respect to which all maps
\[
\Ga\ni\ga\mapsto\skg{f,\ga}:=\int_\X f\sk{x}\ga\sk{dx}=\sum_{x\in\ga}f\sk{x},\quad f\in C_0\sk{\X},
\]
are continuous. Here, $C_0\sk{\X}$ is
the space of all continuous real-valued functions on $\X$
with compact support. We will denote by $\B\sk{\Ga}$ the
Borel $\sigma$-algebra on $\Ga$. We note that $\Ga$ being endowed
with the vague  topology is a Polish space, see e.g.\ \cite{MKM}.

A pair potential is a Borel-measurable function $\phi:\X\goto\R\cup\skf{+\infty}$
such that $\phi\sk{-x}=\phi\sk{x}\in\R$ for all $x\in\X\setminus\skf{0}$.
For $\ga\in\Ga$ and $x\in\X\setminus\ga$, we define a relative
energy of interaction between a particle at $x$ and the
configuration $\ga$ as follows:
\[
E\sk{x,\ga}:=\left\{
\begin{aligned}
&\sum_{y\in\ga}\phi\sk{x-y},&& \text{if } \sum_{y\in\ga}\md{\phi\sk{x-y}}
< +\infty ,\\
& +\infty, &&\text{otherwise.}
\end{aligned}
\right.
\]

A probability measure $\mu$ on $\sk{\Ga,\B\sk{\Ga}}$ is
called a (grand canonical) Gibbs measure corresponding to
the pair potential $\phi$ and activity $z>0$ if it satisfies
the Georgii--Nguyen--Zessin identity (\cite[Theorem 2]{NZ},
see also \cite[Theorem 2.2.4]{Kuna}):
\begin{equation}\label{GNZ}
\int_\Ga\mu\sk{d\ga}\int_\X\ga\sk{dx} F\sk{\ga,x}
=\int_\Ga \mu\sk{d\ga}\int_\X z\,dx \exp\left[-E\sk{x,\ga}\right]
F\sk{\ga\cup x,x}
\end{equation}
for any measurable function $F:\Ga\times\X\goto\left[0;+\infty\right]$.
We denote the set of all such measures $\mu$ by $\G \sk{z,\phi}$.

Let us formulate  conditions on the pair potential $\phi$.

(S) (Stability) There exists $B\geq0$ such that, for any $\ga\in\Ga$, $\md{\ga}<\infty$,
\[
\sum_{\skf{x,y}\subset\ga}\phi\sk{x-y}\geq -B\md{\ga}.
\]

In particular, condition (S) implies that $\phi\sk{x}\geq -2B$,
$x\in\X$.

(P) (Positivity) We have
\[
\phi\sk{x}\geq 0, \quad x\in\X.
\]

The condition (P) is stronger than (S). More precisely, if (P)
holds, then we can choose $B=0$ in (S).

(LA-HT) (Low activity-high temperature regime) We have:
\[
\int_\X\md{e^{-\phi\sk{x}}-1}z\,dx<\sk{2e^{1+2B}}^{-1},
\]
where $B$ is as in (S).

In particular, if (P) holds, then (LA-HT) means:
\[
\int_\X\md{e^{-\phi\sk{x}}-1}z\,dx<\sk{2e}^{-1}.
\]

Let $\mu\in\G \sk{z,\phi}$. Assume that, for any $n\in\N$, there
exists a non-negative, measurable symmetric function
$k_\mu^{\sk{n}}$ on $\sk{\X}^n$ such that, for any measurable
symmetric function $f^{\sk{n}} :\sk{\X}^n \goto
\left[0,+\infty\right]$
$$
\int_\Ga\skg{f^{\sk{n}},:\ga^{\otimes n}:}\,\mu\sk{d\ga}
=\frac{1}{n!}\int_{\sk{\X}^n}f^{\sk{n}}
\sk{x_1,\ldots,x_n}k_\mu^{\sk{n}}
\sk{x_1,\ldots,x_n}\,d x_1 \dotsm d x_n.
$$
Here
\[
\skg{f^{\sk{n}},:\ga^{\otimes n}:}:=
\sum_{\skf{x_1,\ldots,x_n}\subset\ga}f^{\sk{n}}
\sk{x_1,\ldots,x_n}.
\]

The functions $k_\mu^{\sk{n}}$ are called correlation functions
of the measure $\mu$. If there exists a constant $\xi>0$
such that
\begin{equation}\label{RB}
\forall \sk{x_1,\ldots,x_n}\in\sk{\X}^n
:\quad k_\mu^{\sk{n}} \sk{x_1,\ldots,x_n}\leq\xi^n,
\end{equation}
then we say that the correlation functions $k_\mu^{\sk{n}}$
satisfy the Ruelle bound.

Under the conditions (S) and (LA-HT), there exists a Gibbs
measure $\mu\in\G \sk{z,\phi}$ which has correlation
functions satisfying the Ruelle bound, see e.g.\ \cite{Ru69}. This
measure $\mu$ is constructed as a weak limit of finite volume
Gibbs measures  with empty boundary condition, see \cite{Mi67} for details. We will call this measure the Gibbs measure
corresponding to $\sk{z,\phi}$ and the construction with
empty boundary condition. 

In what follows, we will always
assume that (S) and (LA-HT) are satisfied and the Gibbs
measure $\mu$ as discussed above is fixed. We note that,
if the condition (P) is satisfied, then this measure $\mu$
is unique in the set $\G \sk{z,\phi}$, see \cite{Ru69} and
\cite[Theorem~6.2]{Kuna2}. 
We also note  that the relative energy $E\sk{x,\ga}$ is finite
$dx\,\mu\sk{d\ga}$-a.e.\ on $\X\times\Ga$.

Via a recursion formula, one can transform the correlation functions
$k_\mu^{\sk{n}}$ into the Ursell functions $u_\mu^{\sk{n}}$ and vice
versa, see e.g.\ \cite{Ru69}. Their relation is given by
\begin{equation}\label{urs_def}
k_\mu\sk{\eta}=\sum u_\mu\sk{\eta_1}\dotsm u_\mu\sk{\eta_j},\quad
\eta\in\Ga_0,\ \eta\ne\varnothing,
\end{equation}
where
\[
\Ga_0:=\skf{\ga\in\Ga:  \md{\ga}<\infty },
\]
for any $\eta=\skf{x_1,\ldots,x_n}\in\Ga_0$
\[
k_\mu\sk{\eta}:=k_\mu^{\sk{n}}\sk{x_1,\ldots,x_n}, \qquad
u_\mu\sk{\eta}:=u_\mu^{\sk{n}}\sk{x_1,\ldots,x_n} ,
\]
and the summation in (\ref{urs_def}) is over all partitions
of the set $\eta$ into nonempty mutually disjoint subsets
$\eta_1,\dots,\eta_j\subset\eta$ such that
$\eta_1\cup\dotsm\cup\eta_j=\eta$, $j\in\N$. For example,
\begin{align*}
k_\mu^{\sk{1}}\sk{x}&=u_\mu^{\sk{1}}\sk{x},\\
k_\mu^{\sk{2}}\sk{x_1,x_2}&=u_\mu^{\sk{2}}\sk{x_1,x_2}
+u_\mu^{\sk{1}}\sk{x_1}u_\mu^{\sk{1}}\sk{x_2}.
\end{align*}

For our fixed Gibbs measure $\mu$, both the correlation
functions and the Ursell functions of $\mu$ are translation
invariant. In particular, the first correlation function
$k_\mu^{\sk{1}}\sk{\cdot}$ is a constant, which we denote
by $k_\mu^{\sk{1}}$.

Furthermore, for any $n\in\N$,
\begin{equation}\label{Urs_int}
U_\mu^{\sk{n+1}}\in L^1\sk{\sk{\X}^n,dx_1\dotsm dx_n},
\end{equation}
where
\begin{equation}\label{Urs_n+1_def}
U_\mu^{\sk{n+1}}\sk{x_1,\dots, x_n}:=u_\mu^{\sk{n+1}}\sk{x_1,\dots, x_n,0}, \qquad  \sk{x_1,\dots, x_n}\in \sk{\X}^n,
\end{equation}
see \cite{wewewe}.

As a straightforward corollary of the Georgii--Nguyen--Zessin
identity (\ref{GNZ}), we get the following equality:
\begin{align}\notag
&\int_\Ga \mu\sk{d\ga}\int_\X\ga\sk{dx_1}\int_\X\ga\sk{dx_2}
F\sk{\ga,x_1,x_2}\\
&=\int_\Ga \mu\sk{d\ga}\int_\X z\, dx_1\int_\X z\, dx_2\exp\left[
-E\sk{x_1,\ga}-E\sk{x_2,\ga}-\phi\sk{x_1-x_2}
\right]\notag\\
&\hphantom{=\int_\Ga \mu\sk{d\ga}\int_\X z\, dx_1\int_\X z\, dx_2} \ \times F\sk{\ga\cup \skf{x_1, x_2},x_1,x_2}\notag\\
&\quad + \int_\Ga \mu\sk{d\ga}\int_\X z\,dx \exp\left[
-E\sk{x,\ga}\right] F\sk{\ga\cup x,x,x}\label{double_GNZ}
\end{align}
for any measurable function
$F:\Ga\times\X\times\X\goto\left[0,+\infty\right]$.

Let $f:\X\goto\R$ be such that $e^f-1\in L^1\sk{\X,dx}$.
Then, using the representation
\[
e^{\skg{f,\ga}}=1+\sum_{n=1}^\infty\skg{\sk{e^f-1}^{\otimes
n},:\ga^{\otimes n}:},
\]
we get
\begin{multline}\label{exp_cor}
\int_\Ga  e^{\skg{f,\ga}} \mu\sk{d\ga} \\= 1+\sum_{n=1}^\infty
\frac1{n!}\,\int_{\sk{\X}^n} \sk{e^f-1}^{\otimes n} \sk{x_1,\ldots,x_n}
k_\mu^{\sk{n}}\sk{x_1,\dots,x_n}\, dx_1\dotsm dx_n.
\end{multline}
Hence, by using the Ruelle bound, we conclude that
$e^{\skg{f,\cdot}}\in L^1\sk{\Ga,\mu}$.
Furthermore, if  $e^{2f}-1\in L^1(\R^d,dx)$, then $e^{\langle f,\cdot\rangle}\in L^2(\Gamma,\mu)$.

\section{Kawasaki and Glauber dynamics}\label{uitgyuiftu}
We introduce the set $\F C_{\mathrm b}\sk{C_0\sk{\X},\Ga}$ of all
functions of the form
\[
\Ga\ni\ga \mapsto F\sk{\ga}
=g_F\sk{\skg{\varphi_1,\gamma},\ldots,\skg{\varphi_N,\gamma}},
\]
where $N\in\N$, $\varphi_1,\ldots,\varphi_N\in C_0\sk{\X}$,
and $g_F\in C_{\mathrm b}\sk{\X}$, where $C_{\mathrm b}\sk{\X}$ denotes the
set of all continuous bounded functions on $\R^N$.

For each function $F:\Ga\goto\R$, $\ga\in\Ga$,
and $x,y\in\X$, we denote
\begin{align*}
\sk{D_x^{-}F}\sk{\ga}&:=F\sk{\ga\setminus x}-F\sk{\ga},\\
\sk{D_{xy}^{-+}F}\sk{\ga}&
:=F\sk{\ga\setminus x\cup y}-F\sk{\ga}.
\end{align*}

We fix a function
$a:\X\goto\left[ 0, +\infty\right)$
such that $a\sk{-x}=a\sk{x}$, $x\in\X$,
and $a\in L^1\sk{\X,dx}$. We define bilinear forms
\begin{align*}
\EE_\eps\sk{F,G}
&:=\frac{1}{2}\int_\Ga\mu\sk{d\ga}\int_\X\ga\sk{dx}
\int_\X dy  \, a_\varepsilon(x-y) \\
&\qquad \times \exp\left[ -E\sk{y,\ga\setminus x}\right]
\sk{D_{xy}^{-+}F}\sk{\ga}\sk{D_{xy}^{-+}G}\sk{\ga}, \qquad
\eps>0,\\
\EE_0\sk{F,G}&:=\alpha \int_\Ga\mu\sk{d\ga}\int_\X\ga\sk{dx}
\sk{D_{x}^{-}F}\sk{\ga}\sk{D_{x}^{-}G}\sk{\ga},
\end{align*}
where $F,G\in\F C_{\mathrm b}\sk{C_0\sk{\X},\Ga}$, $a_\eps(\cdot)$ is defined by \eqref{rtsra}, and $\alpha$ is given by \eqref{iyuyd}.

The next theorem follows from \cite[Proposition~3.1 and  Theorem~3.1]{KL} and \cite[Proposition~4.3]{KLR}.

\begin{theorem}\label{Hunt_thm}
\rom{i)} For each $\eps\geq 0$, the bilinear form
$\sk{\EE_\eps,\F C_{\mathrm b}\sk{C_0\sk{\X},\Ga}}$
is closable on $L^2\sk{\Ga,\mu}$ and its closure
will be denoted by $\sk{\EE_\eps, \dom\sk{\EE_\eps}}$.

\rom{ii)} Denote by $(H_\eps, \dom(H_\eps))$, $\eps\geq 0$, the generator of
$\sk{\EE_\eps, \dom\sk{\EE_\eps}}$. Then
\[
\F C_{\mathrm b}\sk{C_0\sk{\X},\Ga}\subset
\bigcap_{\eps\geq 0}\dom\sk{H_\eps},
\]
and for any $F\in\F C_{\mathrm b}\sk{C_0\sk{\X},\Ga}$
\begin{align}
\sk{H_\eps F}\sk{\ga}=&-\int_\X\ga\sk{dx}\int_\X dy\,  
a_\varepsilon(x-y)\notag\\
&\qquad\qquad\qquad
\times\exp\left[-E\sk{y,\ga\setminus x}\right]\sk{D_{xy}^{-+}F}\sk{\ga},\quad \eps>0,\label{act_H_Kaw}\\
\sk{H_0 F}\sk{\ga}=&-\alpha\int_\X\ga\sk{dx}\sk{D_x^{-}F}\sk{\ga}
-\alpha\int_\X z \, dx \exp\left[ -E\sk{x,\ga}\right]
\sk{D_x^{+}F}\sk{\ga}.\label{act_H_Glaub}\end{align}

\rom{iii)} For each $\eps\geq 0$, there exists a conservative
Hunt process
\[
\pmb{M}^\eps=\sk{\pmb{\Omega}^\eps, \ \pmb{F}^\eps, \
\sk{\pmb{F}^\eps_t}_{t\geq 0}, \ \sk{\pmb{\Theta}^\eps_t}_{t\geq 0}, \ \sk{\pmb{X}^\eps\sk{t}}_{t\geq 0}, \ \sk{\pmb{P}^\eps_\ga}_{\ga\in\Ga}}
\]
on $\Ga$ \rom(see e.g.  \rom{\cite[p.~92]{MaRo})} which is properly associated
with $\sk{\EE_\eps, \dom\sk{\EE_\eps}}$, i.e., for all \rom($\mu$-versions
of\,\rom) $F\in L^2\sk{\Ga,\mu}$ and all $t>0$ the function
\[
\Ga\ni\ga\mapsto\sk{p_t^\eps F}\sk{\ga}
:=\int_\Omega F\sk{\pmb{X}^\eps\sk{t}}d\pmb{P}^\eps_\ga
\]
is an $\EE_\eps$-quasi-continuous version of
$\exp\left[-tH_\eps\right]F$. $\pmb{M}^\eps$ is up to
$\mu$-equivalence unique \rom(cf.\ \rom{\cite[Chap.\ IV, Sect.\ 6]{MaRo})}.
In particular, $\pmb{M}^\eps$ has $\mu$ as  invariant
measure.
\end{theorem}
\begin{remark}
In Theorem~\ref{Hunt_thm}, $\pmb{M}^\eps$ can be taken canonical,
i.e., $\pmb{\Omega}^\eps$ is the set
$D\sk{\left[0,+\infty\right),\Gamma}$ of all {\it c\'adl\'ag\/} functions
$\omega:\left[0,+\infty\right)\goto\Ga$ (i.e., $\omega$
is right continuous on $\left[0,+\infty\right)$ and has
left limits on $\sk{0,+\infty}$),
$\pmb{X}^\eps\sk{t}\sk{\omega}=\omega\sk{t},t\geq0$,
$\omega\in\pmb{\Omega}^\eps$, $\sk{\pmb{F}^\eps_t}_{t\geq0}$
together with $\pmb{F}^\eps$ is the correponding minimum
completed admissible family (cf.\ \cite[Section 4.1]{Fu80}) and $\pmb{\Theta}^\eps_t$, $t\geq0$, are the corresponding
natural time shifts.
\end{remark}

\section{Convergence of the generators}\label{section4}
We will now study the limiting behavior of the generators
of the Kawasaki dynamics, $H_\eps$, as $\eps\goto
0$. We start with the following
\begin{lemma}\label{exp_in_dom}
For any $\eps\geq 0$ and any $\varphi\in C_0\sk{\X}$, the
function $F\sk{\ga}:=e^\skg{\varphi,\ga}$
belongs to $\dom\sk{H_\eps}$
and the action of $H_\eps$ on $F$ is given by  formula
\eqref{act_H_Kaw} for $\eps>0$ and by \eqref{act_H_Glaub}
for $\eps=0$.
\end{lemma}
\begin{proof}
We first note that since $e^{2\varphi}-1\in L^1\sk{\X,dx}$, we have
$e^{\skg{\varphi,\cdot}}\in L^2 \sk{\Ga,\mu}$.

Assume that $\eps>0$. For each $n\in\N$, we define $g_n\in
C_{\mathrm b}\sk{\R}$ by
\begin{equation}\label{uihgiulgygil}
g_n\sk{u}=\left\{
\begin{aligned}
&e^u, &u\leq n, \\
&e^n, &u>n.
\end{aligned}
\right.
\end{equation}
Then
$
g_n(\langle\varphi,\cdot\rangle)
\in\F C_{\mathrm b}\sk{C_0\sk{\X},\Ga}$.
Since
\[
g_n\sk{\skg{\varphi,\ga}}\leq e^{\skg{\varphi,\ga}}, \quad
\ga\in\Ga,
\]
by the majorized convergence theorem, we have
\[
g_n\sk{\skg{\varphi,\cdot}}\goto e^{\skg{\varphi,\cdot}}
\text{ in } L^2\sk{\Ga,\mu} \text{ as } n\goto\infty.
\]
Next, by  (\ref{double_GNZ}), (\ref{exp_cor}), and the Ruelle bound,
\begin{align*}
&\iGG \left(-\iXg \iXy\, a_\varepsilon (x-y) \right.
\left.\ \exp\skp{-E\sk{y,\ga\setminus x}+\skg{\varphi,\ga}}
\sk{e^{-\varphi\sk{x}+\varphi\sk{y}}-1}\vphantom{\iX}\right)^2\\
&\quad=\iX z\, dx_1 \iX dy_1 \iX z\, dx_2 \iX  dy_2\, \aeps{x_1-y_1}
\aeps{x_2-y_2} \\
&\quad\qquad\times\sk{e^{-\varphi\sk{x_1}}\sk{e^{\varphi\sk{y_1}}-1}
+\sk{e^{-\varphi\sk{x_1}}-1}}
\sk{e^{-\varphi\sk{x_2}}\sk{e^{\varphi\sk{y_2}}-1}
+\sk{e^{-\varphi\sk{x_2}}-1}}\\
&\quad\qquad\times\exp\skp{-\phi\sk{x_1-x_2}-\phi\sk{y_1-x_2}
-\phi\sk{y_2-x_1}}\\
&\quad\qquad\times\iGG\exp\skp{-\sum_{u\in\ga}\sk{\phi\sk{u-x_1}+
\phi\sk{u-x_2}+\phi\sk{u-y_1}+\phi\sk{u-y_2}}}\\
&\quad\quad +\iXx \int_{\R^d}dy_1\int_{\R^d}dy_2
\, \aeps{x-y_1}
\aeps{x-y_2} \\
&\quad\quad\times\sk{e^{-\varphi\sk{x}}\sk{e^{\varphi\sk{y_1}}-1}
+\sk{e^{-\varphi\sk{x}}-1}}\sk{e^{-\varphi\sk{x}}\sk{e^{\varphi\sk{y_2}}-1}
+\sk{e^{-\varphi\sk{x}}-1}}\\
&\quad\quad\times\iGG\exp\skp{-\sum_{u\in\ga}\sk{\phi(u-x)+\phi\sk{u-y_1}+\phi\sk{u-y_2}}}\\
&\quad<\infty.
\end{align*}
Here, we used the following estimate: for any $x,y_1,y_2\in\X$
\begin{align}
&\iX \md{e^{-\phi\sk{u-x}-\phi\sk{u-y_1}-\phi(u-y_2)}-1}du \notag\\
&\qquad=\iX\big| e^{-\phi\sk{u-x}-\phi(u-y_1)}\sk{e^{-\phi\sk{u-y_2}}-1}\notag\\ &\qquad\quad
+e^{-\phi\sk{u-x}}\sk{e^{-\phi\sk{u-y_1}}-1}+\sk{e^{-\phi\sk{u-x}}-1}\big|du\notag\\
&\qquad\leq\sk{e^{4B}+e^{2B}+1}\iX\md{e^{-\phi\sk{u}}-1}du<\infty,\label{spec_est}
\end{align}
where $B$ is as in (S), and an analogous estimate for the
function
\[
e^{-\phi\sk{u-x_1}-
\phi\sk{u-x_2}-\phi\sk{u-y_1}-\phi\sk{u-y_2}}-1.
\]

By using \eqref{uihgiulgygil}, we get, for any $\ga\in\Ga$, $x\in\ga$, and $y\in\X$:
\begin{align*}
&\md{g_n\sk{\skg{\varphi,\ga\setminus x\cup y}}-
g_n\sk{\skg{\varphi,\ga}}}\\
&\qquad=\md{g_n\sk{\skg{\varphi,\ga}-\varphi\sk{x}+\varphi\sk{y}}-
g_n\sk{\skg{\varphi,\ga}}}\\
&\qquad\leq \exp\skp{\max\skf{\skg{\varphi,\ga}-\varphi\sk{x}+\varphi\sk{y},
\skg{\varphi,\ga}}}\sk{-\varphi\sk{x}+\varphi\sk{y}}\\
&\qquad\leq \exp\skp{\skg{\md{\varphi},\ga}+\md{\varphi\sk{x}}+
\md{\varphi\sk{y}}}\sk{|\varphi\sk{x}|+|\varphi\sk{y}|},
\end{align*}
and, hence, for any $\ga\in\Ga$,
\begin{align}
&\md{-\iXg\iXy\,\aeps{x-y}\exp\skp{-E\sk{y,\ga\setminus x}}
D_{xy}^{-+}g_n\sk{\skg{\varphi,\ga}}\vphantom{\iX}}\notag\\
&\qquad\leq \iXg\iXy\,\aeps{x-y}\exp\skp{-E\sk{y,\ga\setminus x}+\skg{\md{\varphi},\ga}+\md{\varphi\sk{x}}+
\md{\varphi\sk{y}}}\notag\\
&\qquad\qquad\times{}\sk{|\varphi\sk{x}|+|\varphi\sk{y}|}.\label{mod_est1}
\end{align}

Analogously to (\ref{spec_est}), we conclude that the right
hand side of (\ref{mod_est1}), as a function of $\ga\in\Ga$,
belongs to $L^2\sk{\Ga, d\mu}$. Therefore, by the majorized
convergence theorem, 
\begin{multline*}
-\iXg\iXy\,\aeps{x-y}\exp\skp{-E\sk{y,\ga\setminus x}}
D_{xy}^{-+}g_n\sk{\skg{\varphi,\ga}}\\
\goto -\iXg\iXy\,\aeps{x-y}\exp\skp{-E\sk{y,\ga\setminus x}} D_{xy}^{-+}e^{\skg{\varphi,\ga}}
\end{multline*}
in $L^2\sk{\Ga,\mu}$ as $n\goto\infty$. From here, the statement
of the lemma follows in the case $\eps>0$. The case $\eps=0$
can be treated analogously.
\end{proof}

We may rewrite $H_\eps=H_\eps^{+} + H_\eps^{-}, \eps\geq 0$, where
\begin{align*}
\sk{H_\eps^{-}F}\sk{\ga}&=-\iXg
\sk{D_x^{-}F}\sk{\ga} \iXy \,  e^{-E\sk{y, \ga \setminus x}}\aeps{x-y}, \\
\sk{H_\eps^{+}F}\sk{\ga}&=- \iXy\iXg  \ e^{-E\sk{y,\ga \setminus x}}\aeps{x-y} \skp{F\sk{ \ga\setminus x\cup y} -F\sk{\ga
\setminus x} } ,\\
\sk{H_0^{-} F}\sk{\ga}&=-\alpha\int_\X\ga\sk{dx}
\sk{D_x^{-}F}\sk{\ga},\\
\sk{H_0^{+} F}\sk{\ga}
&=-\alpha\iXy \, \exp\left[ -E\sk{y,\ga}\right]
\sk{D_y^{+}F}\sk{\ga}.
\end{align*}

\begin{theorem}\label{gen_conv}
Assume that the pair potential $\phi$ and activity $z>0$
satisfy the conditions \rom{(S)} and \rom{(LA-HT)}.  Let
$\mu$ be the Gibbs measure from $\G\sk{z,\phi}$ which corresponds
to the construction with empty boundary condition. Then,
for any $\varphi\in C_0\sk{\X}$,
\[
H_\eps^\pm e^{\skg{\varphi,\cdot}}\goto H_0^\pm e^{\skg{\varphi,\cdot}}
\text{ \rom{in} } L^2\sk{\Ga,\mu} \text{ \rom{as} } \eps\goto 0,
\]
so that 
\[
H_\eps e^{\skg{\varphi,\cdot}}\goto H_0 e^{\skg{\varphi,\cdot}}
\text{ \rom{in} } L^2\sk{\Ga,\mu} \text{ \rom{as} } \eps\goto 0,
\]
\end{theorem}
\noindent{\it Proof}. We fix any $\varphi\in C_0\sk{\X}$ and denote
$F\sk{\ga}:=e^{\skg{\varphi,\ga}}$. We need to prove that 
\begin{align}
\iG\sk{H_\eps^\pm F}^2\sk{\ga}\,\mu\sk{d\ga} & \goto
\iG\sk{H_0^\pm F}^2\sk{\ga}\,\mu\sk{d\ga} \text{ \ as } \eps\goto 0,\label{Apm}\\
\iG\sk{H_\eps^\pm F}\sk{\ga}\sk{H_0^\pm F}\sk{\ga}\,\mu\sk{d\ga} & \goto
\iG\sk{H_0^\pm F}^2\sk{\ga}\,\mu\sk{d\ga} \text{ \ as } \eps\goto 0.\label{Bpm}
\end{align}

Using (\ref{GNZ}) and (\ref{double_GNZ}), we get
\begin{align}&\iG\sk{H_0^{-}F}^2\sk{\ga}\,\mu\sk{d\ga}\notag\\
&\qquad=\alpha^2\iXx\f{x}^2\iGG\exp\skp{-E\sk{x,\ga}+\skg{2\varphi,\ga}}\notag\\
&\qquad\quad
+\alpha^2 \iXx\iXx'\,\f{x}\g{x}\f{x'}\g{x'}e^{-\phi\sk{x-x'}}\notag \\
&\qquad\qquad\times\iGG\exp\skp{-E\sk{x,\ga}-E\sk{x',\ga}+\skg{2\varphi,\ga}}
\label{m-fin}
\end{align}
and
\begin{align}
&\iG\sk{H_0^{+}F}^2\sk{\ga}\,\mu\sk{d\ga}\notag\\
&\qquad=\alpha^2 \iXx\iXx'\,\f{x}\f{x'}\notag\\
&\qquad\qquad\times\iGG\exp\skp{-E\sk{x,\ga}-E\sk{x',\ga}+\skg{2\varphi,\ga}}.
\label{p-fin}
\end{align}

Using the same arguments and changing the variables, we have
\begin{align}
&\iG\sk{H_\eps^{-}F}^2\sk{\ga}\,\mu\sk{d\ga}\notag\\
&\qquad=\iXx\iXy \iXy'\,\f{x}^2 \aeps{x-y} \aeps{x-y'}\notag \\
&\qquad\qquad\times \iGG  \,
 \exp\skp{-E\sk{x,\ga}-E\sk{y, \ga}-E\sk{y', \ga }+\skg{2\varphi,\ga}}\notag\\
&\qquad\quad+ \iXx\iXx'\iXy \iXy'\,\f{x}\g{x}\f{x'}\g{x'}\notag\\
&\qquad\qquad\times \aeps{x-y} \aeps{x'-y'} e^{-\phi\sk{x-x'}}e^{-\phi\sk{x-y'}}e^{-\phi\sk{y-x'}}\notag\\
&\qquad\qquad\times\iGG\exp\skp{-E\sk{x,\ga}-E\sk{x',\ga}-E\sk{y, \ga}-E\sk{y', \ga }+\skg{2\varphi,\ga}}\notag\\
&\qquad=\iXx\iXy \iXy'\,\f{x}^2 a\sk{y} a\sk{y'}\notag \\
&\qquad\qquad\times \iGG  \,
 \exp\skp{-E\sk{x,\ga}-E\sk{\frac{y}{\eps}+x, \ga}
 -E\sk{\frac{y'}{\eps}+x, \ga }+\skg{2\varphi,\ga}}\notag\\
&\qquad\quad+ \iXx\iXx'\iXy \iXy'\,\f{x}\g{x}\f{x'}\g{x'}\notag\\
&\qquad\qquad\times a\sk{y} a\sk{y'} e^{-\phi\sk{x-x'}}
e^{-\phi\sk{\frac{y'}{\eps}+x'-x}}e^{-\phi\sk{\frac{y}{\eps}+x-x'}}\notag\\
&\qquad\quad\times\iGG\exp\bigg[-E\sk{x,\ga}-E\sk{x',\ga}\notag\\
&\qquad\qquad
-E\sk{\frac{y}{\eps}+x, \ga}-E\sk{\frac{y'}{\eps}+x', \ga }+\skg{2\varphi,\ga}\bigg].\label{aaa}
\end{align}
Analogously,
\begin{align}
&\iG\sk{H_\eps^{+}F}^2\sk{\ga}\,\mu\sk{d\ga}\notag\\
&\qquad=\iXx\iXy \iXy'\,\f{x}\f{\frac{y'-y}{\eps}+x} a\sk{y} a\sk{y'}\notag \\
&\quad\quad\times \iGG  \,
 \exp\skp{-E\sk{x-\frac{y}{\eps},\ga}-E\sk{x, \ga}
 -E\sk{\frac{y'-y}{\eps}+x, \ga }+\skg{2\varphi,\ga}}\notag\\
&\qquad+ \iXx\iXx'\iXy \iXy'\,\f{y}\g{\frac{x}{\eps}+y}\f{y'}\g{\frac{x'}{\eps}+y'}\notag\\
&\quad\qquad\times a\sk{x} a\sk{x'} e^{-\phi\sk{\frac{x-x'}{\eps}+y-y'}}
e^{-\phi\sk{\frac{x}{\eps}+y-y'}}e^{-\phi\sk{y-\frac{x'}{\eps}-y'}}\notag\\
&\quad\qquad\times\iGG\exp\bigg[-E\sk{\frac{x}{\eps}+y,\ga}\notag\\
&\qquad\qquad-E\sk{\frac{x'}{\eps}+y',\ga}-E\sk{y, \ga}-E\sk{y', \ga }+\skg{2\varphi,\ga}\bigg].\label{bbb}
\end{align}
In the same way,
\begin{align}
&\iG\sk{H_\eps^{-}F}\sk{\ga}\sk{H_0^{-}F}\sk{\ga}\,\mu\sk{d\ga}\notag\\
&\quad=\alpha\iXx\iXy\f{x}^2 a\sk{y} \notag \\
&\quad\qquad\times \iGG  \,
 \exp\skp{-E\sk{x,\ga}-E\sk{\frac{y}{\eps}+x, \ga}+\skg{2\varphi,\ga}}\notag\\
&\quad\quad+\alpha\ \iXx\iXx'\iXy\f{x}\g{x}\f{x'}\g{x'}\notag\\
&\quad\qquad\qquad\times a\sk{y}  e^{-\phi\sk{x-x'}}e^{-\phi\sk{\frac{y}{\eps}+x-x'}}\notag\\
&\quad\qquad\times\iGG\exp\skp{-E\sk{x,\ga}-E\sk{x',\ga}-E\sk{\frac{y}{\eps}+x, \ga}+\skg{2\varphi,\ga}},\label{ccc}
\end{align}
and finally,
\begin{align}
&\iG\sk{H_\eps^{+}F}\sk{\ga}\sk{H_0^{+}F}\sk{\ga}\,\mu\sk{d\ga}\notag\\
&\qquad=\alpha\iXx\iXy \iXy'\,\f{y}\f{y'} \g{\frac{x}{\eps}+y} e^{-\phi\sk{\frac{x}{\eps}+y-y'}} a\sk{x} \notag\\
&\qquad\qquad\times  \iGG  \,
 \exp\skp{-E\sk{\frac{x}{\eps}+y,\ga}-E\sk{y, \ga}-E\sk{y', \ga }+\skg{2\varphi,\ga}}.\label{ddd}
\end{align}

Since the functions $e^\varphi-1, e^{-\phi}-1$ are bounded and
integrable on $\X$, from the Ruelle bound  and \eqref{exp_cor} we conclude that all the 
integrals over $\Ga$ in the right hand sides of the equalities \eqref{aaa}--\eqref{bbb}
are bounded by constants. Therefore, by the majorized convergence theorem, to find the limit of these expressions as $\varepsilon\to0$, it suffices to find the limit of the corresponding integrals over $\Gamma$ for fixed variables $x$, $x'$, $y$, $y'$.    
Therefore, due to \eqref{m-fin} and \eqref{p-fin}, formulas
\eqref{Apm}, \eqref{Bpm} will immediately follow from the
following

\begin{lemma}\label{mainlemma} Let a function $\psi:\R^d\to \R$ be such that $e^\psi-1$ is bounded and integrable.  Suppose that $x,y,x',y'\in\R^d$ and  $x\ne y$.  Then
\begin{align}
\int_{\Gamma}  \exp\skp{-E\sk{\frac{x}{\eps}+x', \ga}
 -E\sk{\frac{y}{\eps}+y', \ga }+\skg{\psi,\ga}}\,\mu(d\gamma)
& \goto \frac{\sk{k^\sk{1}_\mu}^2}{z^2} \int_\Gamma
 \exp\skp{\skg{\psi,\ga}}\,\mu(d\gamma),\label{4a}\\
 \int_\Gamma  \exp\skp{-E\sk{\frac{x}\eps+x', \ga}
 +\skg{\psi,\ga}}\,\mu(d\gamma)
& \goto \frac{k^\sk{1}_\mu}{z} \iG
 \exp\skp{\skg{\psi,\ga}}\,\mu(d\gamma)\label{4b}
\end{align}
as $\eps\to0$.
\end{lemma}

\begin{proof}
By \eqref{exp_cor}, 
\begin{align}
&\int_{\Gamma}  \exp\skp{-E\sk{\frac{x}{\eps}+x', \ga}
 -E\sk{\frac{y}{\eps}+y', \ga }+\skg{\psi,\ga}}\,\mu(d\gamma)
\notag\\
&=1+\sum_{n=1}^\infty\frac{1}{n!}\int_{\sk{\X}^n} \left( \exp\left[
-\phi\sk{\cdot-x(\eps)}-\phi\sk{\cdot-y(\eps)} +\psi\sk{\cdot}
\right]-1\right)^{\otimes n}\sk{u_1,\ldots,u_n} \notag \\
&\qquad \qquad \qquad \times k_\mu^\sk{n} \sk{u_1,\ldots,u_n}\, d
u_1\dotsm du_n, \label{dop}
\end{align}
where $x(\eps):=\frac{x}\eps+x'$, $y(\eps):=\frac{y}\eps+y'$. 

Using the Ruelle bound, (S), and (LA-HT),  we
conclude that, in order to find the limit of the right hand side of
\eqref{dop} as $\eps\goto 0$, it suffices to find the limit of each
term
\begin{align}
C_\eps^\sk{n}:&=\int_{\sk{\X}^n} \left(\exp\left[
-\phi\sk{\cdot-x\sk{\eps}} -\phi\sk{\cdot-y\sk{\eps}}
+\psi\sk{\cdot} \right]-1\right)^{\otimes
n}\sk{u_1,\ldots,u_n}\notag\\
&\qquad\qquad\qquad\qquad\qquad\qquad\qquad\times k_\mu^\sk{n}
\sk{u_1,\ldots,u_n}\, d u_1\dotsm du_n \notag\\
&=\sum_{n_1+n_2+n_3=n}\binom{n}{n_1\,n_2\,n_3}\int_{\sk{\X}^n}
\sk{f_{1,\eps}^{\otimes n_1}\otimes f_{2,\eps}^{\otimes n_2}\otimes f_{3,\eps}^{\otimes n_3}}\sk{u_1,\ldots,u_n}\notag\\
&\qquad\qquad\qquad\qquad\qquad\qquad\qquad\times k_\mu^\sk{n}
\sk{u_1,\ldots,u_n} \,d u_1\dotsm du_n,\label{C_n_eps}
\end{align}
where
\begin{align*}
f_{1,\eps}\sk{u}&:=\sk{\exp\skp{-\psi\sk{u}}-1}
\exp\skp{-\phi\sk{u-x\sk{\eps}}-\phi\sk{u-y\sk{\eps}}},\\
f_{2,\eps}\sk{u}&:=\sk{\exp\skp{-\phi\sk{u-x\sk{\eps}}}-1}
\exp\skp{-\phi\sk{u-y\sk{\eps}}},\\
f_{3,\eps}\sk{u}&:=\exp\skp{-\phi\sk{u-y\sk{\eps}}}-1, \qquad
u\in\R.
\end{align*}

Using \eqref{urs_def}, we see that
\begin{multline*}
\iXn{n}\sk{f_{1,\eps}^{\otimes n_1}\otimes f_{2,\eps}^{\otimes n_2}\otimes f_{3,\eps}^{\otimes n_3}}\vect{u}{n}
k_\mu^\sk{n}\vect{u}{n} \dvect{u}{n} \\
= \sum \iXn{n}\sk{f_{1,\eps}^{\otimes n_1}\otimes
f_{2,\eps}^{\otimes n_2}\otimes f_{3,\eps}^{\otimes n_3}}\vect{u}{n}
u_\mu\sk{\eta_1}\dotsm u_\mu\sk{\eta_j} \dvect{u}{n},
\end{multline*}
where the summation is over all partitions of
$\eta=\skf{u_1,\ldots,u_n}$.
We now have to distinguish the three following cases.

Case 1: Each element $\eta_i$ of the partition is either a subset of
$\skf{u_1,\ldots,u_{n_1}}$, or a subset of
$\skf{u_{n_1+1},\ldots,u_{n_1+n_2}}$, or a subset of
$\skf{u_{n_1+n_2+1},\ldots,u_{n}}$. 

By using the translation
invariance of the Ursell functions, we get that the corresponding
term is equal to
\begin{equation}\label{term1}
\iXn{n}\sk{f_{1,\eps}^{\otimes n_1}\otimes g_{2,\eps}^{\otimes n_2}\otimes g_{3,\eps}^{\otimes n_3}}\,\vect{u}{n}
u_\mu\sk{\eta_1}\dotsm u_\mu\sk{\eta_j} \dvect{u}{n},
\end{equation}
where
\begin{align*}
g_{2,\eps}\sk{u}&:=\sk{\exp\skp{-\phi\sk{u}}-1}
\exp\skp{-\phi\sk{u+\frac{x-y}{\eps}+x'-y'}},\\
g_{3,\eps}\sk{u}&:=\exp\skp{-\phi\sk{u}}-1, \qquad u\in\R.
\end{align*}
Since $e^{-\phi}-1\in L^1(\R^d,dx)$, for each $\delta>0$
there exists $r>0$ such that \begin{equation}\label{yudeyed}\int_{\mathbb R^d} I(|x|\ge r,\, |\phi(x)|\ge\delta)\,dx<\delta,
\end{equation}
where $I(\cdot)$ denotes the indicator function. Hence, by the majorized convergence theorem, 
since $x\ne y$, \eqref{term1} converges
to
\begin{multline*}
\iXn{n} \sk{\exp\skp{\psi\sk{\cdot}}-1}^{\otimes n_1} \\ \otimes
\sk{\exp\skp{-\phi\sk{\cdot}}-1}^{\otimes\sk{n_2+n_3}} \vect{u}{n}
u_\mu\sk{\eta_1}\ldots u_\mu\sk{\eta_j} \dvect{u}{n}.
\end{multline*}

Case 2: There is an element of the partition which has non-empty
intersections with both sets $\skf{u_1,\ldots,u_{n_1}}$
and $\skf{u_{n_1+1},\ldots,u_n}$. 

By  using \eqref{Urs_int},
 \eqref{Urs_n+1_def}, \eqref{yudeyed}, and the majorized convergence theorem,
we conclude that the term converges to zero as $\eps\goto0$.

Case 3: Case 2 is not satisfied, but there is an element $\eta_l$ of
the partition which has non-empty intersections with both sets
$\skf{u_{n_1+1},\ldots,u_{n_1+n_2}}$ and
$\skf{u_{n_1+n_2+1},\ldots,u_{n}}$. 

Shift all the variables
entering $\eta_l$ by $y\sk{\eps}$. Now, analogously to Case~2, the
term converges to zero as $\eps\goto 0$.

Thus, again using \eqref{urs_def}, for each $n\in\N$,
\begin{align*}
C_\eps^\sk{n}&\goto \sum_{n_1+n_2+n_3=n}\binom{n}{n_1\,n_2\,n_3}
\\&\qquad\times\iXn{n_1}\sk{\exp\skp{\psi\sk{\cdot}}-1}
^{\otimes n_1} \vect{u}{n_1}\\&\qquad\qquad\qquad\times k_\mu^\sk{n_1}\vect{u}{n_1} \dvect{u}{n_1}\\
&\qquad\times
\iXn{n_2}\sk{\exp\skp{-\phi\sk{\cdot}}-1}^{\otimes n_2}
\sk{u_{n_1+1},\ldots,u_{n_1+n_2}}\\
&\qquad\qquad\qquad\times k_\mu^\sk{n_2}
\sk{u_{n_1+1},\ldots,u_{n_1+n_2}}\,d u_{n_1+1} \dotsm du_{n_1+n_2}\\
&\qquad\times \iXn{n_3}
\sk{\exp\skp{-\phi\sk{\cdot}}-1}^{\otimes n_3}
\sk{u_{n_1+n_2+1},\ldots,u_{n}}\\
&\qquad\qquad\qquad\times k_\mu^\sk{n_3}
\sk{u_{n_1+n_2+1},\ldots,u_{n}}\, d u_{n_1+n_2+1} \dotsm du_{n}.
\end{align*}
Therefore, the right hand side of \eqref{dop} converges to
\[
\sk{\int_\Ga\exp\skp{-\sum_{u\in\ga}\phi\sk{u}}\,\mu\sk{d\ga}}^2
 \iGG\exp\skp{\skg{\psi,\ga}}.
\]

Let $f\in C_0\sk{\X}$, $f\geq 0$, $f\neq 0$. Then
\begin{align*}
k_\mu^\sk{1}\iX f\sk{x}\,dx&=\iGG\iXg f\sk{x}\\
&=\iGG\iXx\exp\skp{-\sum_{u\in\ga}\phi\sk{u-x}}f\sk{x}\\
&=\iXx \, f\sk{x}\iGG \exp\skp{-\sum_{u\in\ga}\phi\sk{u-x}}\\
&=\iXx\, f\sk{x}\iGG \exp\skp{-\sum_{u\in\ga}\phi\sk{u}}.
\end{align*}
Hence,
\[
\int_\Ga\exp\skp{-\sum_{u\in\ga}\phi\sk{u}}\,\mu\sk{d\ga}=\frac{k_\mu^\sk{1}}{z},
\]
which proves  \eqref{4a}. The proof of \eqref{4b} is analogous. \end{proof}


\section{Convergence of the processes}\label{section5}
For each $\eps\geq 0$, we take the canonical version of the
process $\pmb{M}^\eps$ from Theorem~\ref{Hunt_thm}, and
define a stochastic process
$\pmb{Y}^\eps=\sk{\pmb Y^\eps\sk{t}}_{t\geq 0}$
whose law is the probability measure on
$D\sk{\left[0,+\infty\right),\Ga}$
given by
\[
\pmb{Q}^\eps:=\iG\pmb{P}^\eps_\ga \, \mu\sk{d\ga}.
\]
By virtue of Theorem~\ref{Hunt_thm}, the process $\pmb{Y}^\eps$
has $\mu$ as  invariant measure.

\begin{theorem}\label{jlhgiu}
Let \rom{(P) and  (LA-HT),}  be satisfied. Then
the finite-dimensional distributions of the process
$\pmb{Y}^\eps$ weakly converge to the finite-dimensional
distributions of $\pmb{Y}^0$ as $\eps\goto 0$.
\end{theorem}

\begin{proof} Fix any $0\le t_1<t_2<\dots<t_n$, $n\in\mathbb N$. For $\eps\ge0$, denote by $\mu_{t_1,\dots,t_n}^\eps$ the 
finite-dimensional  distribution of the process  $\pmb{Y}^\eps$ at times $t_1,\dots,t_n$, which is a probability measure on $\Gamma^n$. 
Since $\Gamma$ is a Polish space, by \cite[Chapter II, Theorem~3.2]{Par}, the  measure $\mu$ is tight on $\Gamma$. Since all the marginal distributions of the measure $\mu_{t_1,\dots,t_n}^\eps$ are $\mu$, we therefore conclude that the set $\{\mu_{t_1,\dots,t_n}^\eps\mid\eps>0\}$ is pre-compact in the space $\mathcal M(\Gamma^n)$ of the probability measures on $\Gamma^n$ with respect to the weak topology, see e.g.\ \cite[Chapter II, Section~6]{Par}. Hence, by \cite[Chapter~3, Theorem~3.17]{Dav},
the statement of the theorem will follow from 
Theorem~\ref{gen_conv} if we show that the set of all finite
linear combinations of the functions of the form
$e^\skg{\varphi,\cdot}$, $\varphi\in C_0\sk{\R}$, constitutes a
core of $\sk{H_0, \dom\sk{H_0}}$.

Under the assumptions of the~theorem, it follows from the
proof of \cite[Theorem~4.1]{KL}  that the set of all finite
sums of the functions of the form $\prod_{i=1}^n \skg{f_i,\cdot}$,
$f_i\in C_0\sk{\X}$, $i=1,\ldots,n$, and constants
forms a core of $\sk{H_0, \dom\sk{H_0}}$. Therefore, by the  
polarization identity (see e.g.\  \cite[Chapter~2, formula~(2.7)]{BK}) the set of all finite linear combinations of the functions
of the form $\skg{f,\cdot}^n$, $f\in C_0\sk{\X}$, $n=0, 1, 2, \dots$, forms a core of $\sk{H_0, \dom\sk{H_0}}$.
Hence, to prove the~theorem, it suffices to show that, for
each $n\in\N$ and $f\in C_0\sk{\X}$, the function $\skg{f,\cdot}^n$
can be approximated by finite linear combinations of functions
$e^\skg{\varphi,\cdot}$, $\varphi\in C_0\sk{\X}$, in the
graph norm of the operator $\sk{H_0, \dom\sk{H_0}}$. This
statement will follows from the two following  lemmas.

\begin{lemma}
For any $\varphi,\psi\in C_0\sk{\X}$ and $n\in\N$,
$\skg{\psi,\cdot}^n e^\skg{\varphi,\cdot}\in\dom\sk{H_0}$.
\end{lemma}

\begin{proof}
The proof of this lemma is absolutely analogous to the proof
of~Lemma~\ref{exp_in_dom}.
\end{proof}

\begin{lemma}\label{ytee}
For any $\varphi\in C_0\sk{\X}$, $t\in\R$, and $n=0,1,2,\ldots$,
\begin{equation}\label{dop_conv}
\skg{\varphi,\cdot}^n \frac{1}{u}
\sk{e^{\sk{t+u}\skg{\varphi,\cdot}}-e^{t\skg{\varphi,\cdot}}}\goto
\skg{\varphi,\cdot}^{n+1} e^{t\skg{\varphi,\cdot}} \text{
\rom{as} } u\goto 0,
\end{equation}
where convergence is in the sense of the graph norm of the
operator $\sk{H_0, \dom\sk{H_0}}$.
\end{lemma}
\begin{proof}
Using the estimate
\[
\md{\frac{1}{u}\sk{e^{u\skg{\varphi,\ga}}-1}}
\leq e^\skg{\md{\varphi},\ga} \skg{\md{\varphi},\ga}, \quad
\md{u}\leq 1, \ \varphi\in C_0\sk{\X},
\]
and the majorized convergence theorem, we easily get the
convergence \eqref{dop_conv} in $L^2\sk{\Ga,\mu}$.

Next, we have the estimate, for $\md{u}\leq 1$
\begin{align*}
&\md{
\iXg\sk{D_x^{-}\skg{\varphi,\cdot}^n
\frac{1}{u}\sk{e^{\sk{t+u}\skg{\varphi,\cdot}}-
e^{t\skg{\varphi,\cdot}}}}\sk{\ga}}\\
&\qquad=\left|\iXg\left[
            \sk{
            \skg{\varphi,\ga\setminus x}^n
            -\skg{\varphi,\ga}^n}
            \frac{1}{u}\sk{e^{\sk{t+u}
            \skg{\varphi,\ga\setminus x}}
             -e^{t\skg{\varphi,\ga\setminus x}}}
            \right.\right.\\
&\qquad\qquad\qquad\qquad\ \ + \skg{\varphi,\ga}^n
\sk{e^{t\skg{\varphi,\ga\setminus x}}
-e^{t\skg{\varphi,\ga}}}
\frac{1}{u}\sk{e^{u\skg{\varphi,\ga\setminus x}}-1}\\
&\qquad\left.\left.\qquad\qquad\qquad\ \ + \skg{\varphi,\ga}^n e^{t\skg{\varphi,\ga}}
\frac{1}{u}\sk{e^{u\skg{\varphi,\ga\setminus x}}-e^{u\skg{\varphi,\ga}}}
\right]\right|\\
&\qquad\leq \iXg\left[
\sum_{k=0}^{n-1} \skg{\md{\varphi},\ga}^k \md{\varphi\sk{x}}^{n-k}
e^{\sk{\md{t}+1}\skg{\md{\varphi},\ga}}\skg{\md{\varphi},\ga}\right.\\
&\qquad\qquad\qquad\qquad\qquad + \skg{\md{\varphi},\ga}^n
e^{t\skg{\md{\varphi},\ga}}t \md{\varphi\sk{x}}
e^{\skg{\md{\varphi},\ga}}\skg{\md{\varphi},\ga}\\
&\qquad\qquad\qquad\qquad\qquad +\skg{\md{\varphi},\ga}^n
e^{t\skg{\md{\varphi},\ga}} e^{\skg{\md{\varphi},\ga}}
\md{\varphi\sk{x}}\left.\vphantom{\sum_{k=0}^{n-1}}\right].
\end{align*}
Therefore, by the majorized convergence theorem
\begin{align*}
&\iXg\sk{D_x^{-}\skg{\varphi,\cdot}^n
\frac{1}{u}\sk{e^{\sk{t+u}\skg{\varphi,\cdot}}-
e^{t\skg{\varphi,\cdot}}}}\sk{\ga}\\
&\qquad\goto \iXg \sk{D_x^{-}\skg{\varphi,\cdot}^{n+1}
e^{t\skg{\varphi,\cdot}}}\sk{\ga}
\quad \text{as } u\goto 0
\end{align*}
in $L^2\sk{\Ga,\mu}$.

Finally, noticing that $\exp\skp{-E\sk{x,\ga}}\leq 1$ due
to (P), we conclude, analogously to the above, that
\begin{align*}
&\iXx \exp\skp{-E\sk{x,\ga}}
\sk{D_x^{+}\skg{\varphi,\cdot}^n
\frac{1}{u}\sk{e^{\sk{t+u}\skg{\varphi,\cdot}}-
e^{t\skg{\varphi,\cdot}}}}\sk{\ga}\\
&\qquad\goto \iXx \exp\skp{-E\sk{x,\ga}}
\sk{D_x^{+}\skg{\varphi,\cdot}^{n+1}
e^{t\skg{\varphi,\cdot}}}\sk{\ga}
\quad \text{as } u\goto 0
\end{align*}
in $L^2\sk{\Ga,\mu}$. From here, the statement of the lemma
follows.
\end{proof}

Using the induction in $n=0,1,2,\dots$, we conclude from Lemma~\ref{ytee}
that any function of the form $\skg{f,\cdot}^{n+1}
e^{t\skg{f,\cdot}}$, $f\in C_0(\R^d)$, may be approximated by finite
linear  combinations of the functions $e^{\langle\varphi,\cdot\rangle}$, $\varphi\in C_0(\R^d)$, 
in the graph norm of $\sk{H_0, \dom\sk{H_0}}$.
Letting $t=0$, we get the needed statement.
\end{proof}

\section{Concluding remarks and open problems}\label{section6}

It is known (cf.\ \cite{rebenko}) that any Gibbs measure of  Ruelle type,
corresponding to a superstable potential  $\phi$ (see \cite{Ru70}), satisfies the generalized Ruelle bound:
\begin{equation}\label{GRB} k_\mu(x_1,\dots,x_n)\le C^n 
\exp\left[-\sum_{1\le i<j\le n}\phi(x_i-x_j)\right]\end{equation} 
for some $C>0$ independent of $n\in\mathbb N$ (compare with the usual Ruelle bound \eqref{RB}). 
Using \eqref{GRB} and harmonic analysis on the configuration space (cf.\ \cite{KK}), it is possible to derive the convergence of the generators, as in Theorem \ref{gen_conv},  under the following assumption on a decay of correlations of the measure $\mu$: 
 For each $n\in\mathbb N$ and for $dx_1\dotsm dx_{n}$-a.e.\ $(x_1,\dots,x_{n})\in\left(\mathbb R^d\right)^{n}$,
\begin{align}u_\mu^{(n+1)}\left(x_1,\dots,x_n,\frac{y}{\eps}\right)&\to0,
\notag\\ 
u_\mu^{(n+2)}\left(x_1,\dots,x_n,\frac{y}{\eps},\frac{y'}{\eps}\right)&\to0\quad 
\text{as }\eps\to0,\label{piooy}\end{align}
where the convergence is in the $dy\,dy'$-measure on each compact set in $\left(\mathbb R^d\right)^2$. By \eqref{Urs_int},
\eqref{Urs_n+1_def}, this assumption is satisfied in the low activity-high temperature regime. It is still an open problem
whether other Gibbs measures of  Ruelle type satisfy \eqref{piooy}. We believe that \eqref{piooy} indeed holds for any Gibbs measure of Ruelle type which is a pure phase. Note that, if this were so, we would derive  the convergence of the equilibrium processes, as in Theorem \ref{jlhgiu}, assuming additionally (P).

Next, let us consider the following generalization of the Kawasaki and Glauber dynamics.
 For any $s\in[0,1]$, let us define bilinear forms
\begin{align}
\EE_\eps^s\sk{F,G}
&:=\frac{1}{2}\int_\Ga\mu\sk{d\ga}\int_\X\ga\sk{dx}
\int_\X dy \, c_\eps^s\sk{x,y,\ga\setminus x}
\sk{D_{xy}^{-+}F}\sk{\ga}\sk{D_{xy}^{-+}G}\sk{\ga}, \
\ \eps>0, \label{jhvhgjvf}\\
\EE_0^s\sk{F,G}&:=\alpha \int_\Ga\mu\sk{d\ga}\int_\X\ga\sk{dx}
c_0^s\sk{x,\ga\setminus x}\sk{D_{x}^{-}F}\sk{\ga}\sk{D_{x}^{-}G}\sk{\ga},\notag
\end{align}
where $F,G\in\F C_{\mathrm b}\sk{C_0\sk{\X},\Ga}$
and 
\begin{align*}
c_\eps^s\sk{x,y,\ga}&:=a_\eps\sk{x-y}\exp \left[ sE\sk{ x,\ga} -\sk{
1-s} E\sk{ y,\ga} \right] , \\
c_0^s\sk{x,\ga}&:=\exp\skp{sE\sk{x,\ga}}.
\end{align*}
In particular, for $s=0$, $\EE_\eps^0=\EE_\eps$ and $\EE_0^0=\EE_0$. 
By \cite{KLR}, each of these bilinear forms leads to an equilibrium Markov processes on $\Gamma$,
which is a Kawasaki dynamics for $\eps>0$, and a Glauber dynamics for $\eps=0$.

Under
the same assumptions as in Theorem~\ref{gen_conv}, it
can proved that, for any $F\sk{\ga}:=e^{\skg{\varphi,\ga}}$
with $\varphi\in C_0\sk{\X}$,
\[
\EE_\eps^s\sk{F,F}\goto\EE_0^s\sk{F,F}\quad \text{as } \eps\goto 0.
\]
The idea of proof is the same as before, we only use
the second statement of Lemma~\ref{mainlemma}.
Moreover, we expect  that  (under an additional assumption if $s\in(1/2,1]$) an analog of
Theorem~\ref{gen_conv} is also true in this case. However, to derive from here the convergence of the corresponding processes, as in Theorem \ref{jlhgiu}, we are still missing a theorem on a core for the 
generator of the closure of the bilinear form $\mathcal E_0^s$.

As for a non-equilibrium dynamics, let us
first mention that there are two possible ways of its understanding: 
One way is to consider a dynamics which has some given initial distribution. Such an understanding of stochastic dynamics is natural for infinite particle systems, as  has been observed in many models 
of mathematical physics. The other way is to consider a dynamics  starting from a given configuration. The latter corresponds to the canonical approach in the theory of Markov processes.

In the case of no interaction between particles, $\phi=0$, it is  possible to prove the `Kawasaki to Glauber' convergence for a non-equilibrium dynamics whose initial distribution has Ursell functions decaying at infinity, see \cite{KLR2} for details. However, we do not expect any type of convergence even of a free dynamics if the dynamics starts
from a given configuration.

Another version of our results may be applied to the
dynamics considered in \cite[Section 5]{BCPP}.  There, the `Kawasaki' dynamics corresponding to
the following bilinear form was studied:
\begin{equation}\label{kjugyufguf}
\EE^N_\La\sk{F,G}
:=\frac{1}{2|\Lambda|}\int_{\Ga_\La^N}\mu^N_\La\sk{d\ga}\int_\La\ga\sk{dx}
\int_\La  dy  \sk{D_{xy}^{-+}F}\sk{\ga}\sk{D_{xy}^{-+}G}\sk{\ga}.
\end{equation}
Here, $\La$ is a measurable bounded domain in $\X$, by $\md{\La}$ we
denote the volume of $\La$, $\Ga_\La^N\subset\Ga$ is the set of all
$N$-point subsets of $\La$ and $\mu^N_\La$ is the canonical Gibbs
measure on $\Ga_\La^N$ corresponding to the potential $\phi$ (see \cite{BCPP} for detail, and note that we have chosen empty boundary condition). 

Now, let $\La\nearrow\X$, $N\goto\infty$ and $\frac{N}{\md{\La}}\goto\rho=
\mathrm{const}$ (the so-called $N/V$ limit). Then, by \cite{GKL}, the measures $\mu_\Lambda^N$ have a limiting point in the weak topology of probability measures on $\Gamma$. This limiting point is a (grand canonical) Gibbs measure $\mu$ corresponding to the potential $\phi$. A heuristic calculation shows that the Kawasaki dynamics corresponding to \eqref{kjugyufguf} converges to the Glauber dynamics with equilibrium measure $\mu$. In this way, we obtain an $N/V$-approximation of the Glauber dynamics on $\Gamma$. 

It was shown in \cite{BCPP} that, under (P),
the generator of \eqref{kjugyufguf} has a spectral gap  which is $\ge$
$$ 1-\frac{3(N-1)}{|\Lambda|}\int_{\R^d}(1-e^{-\phi(x)})\,dx,$$
provided that the above value is positive. Hence, we should expect that the generator of the limiting Glauber dynamics has a spectral gap which is $\ge$ 
$$ 1-3\rho\int_{\R^d}(1-e^{-\phi(x)})\,dx. $$
This can be compared with the lower bound of the spectral gap
$$ 1-z\int_{\R^d}(1-e^{-\phi(x)})\,dx, $$
obtained in \cite{KL}.

There is still an open problem whether an equilibrium Glauber dynamics in infinite volume can have a spectral gap if the pair potential $\phi$ has a negative part. We hope that the finite-volume approximations as discussed above should give an insight into this problem. 

It should also be possible to approximate the Glauber dynamics in infinite volume by Glauber dynamics in a finite volume which have a grand canonical Gibbs measure as invariant measure.

There is another interesting open problem in this direction: approximation of the Kawasaki dynamics in infinite volume by finite-volume Kawasaki dynamics. Consider e.g.\ the bilinear form $\mathcal E_1$ as in Section \ref{uitgyuiftu}.  This form can be approximated by the forms \[
\EE_\La\sk{F,F}=\frac{1}{2}\int_\GL\mu_\La\sk{d\ga}\int_\La\ga\sk{dx}
\int_\La dy \ a(x-y) \exp\left[ -E\sk{y,\ga\setminus x}\right]
\md{\sk{D_{xy}^{-+}F}\sk{\ga}}^2.
\]
Here, $\Gamma_\Lambda$ is the set of all finite subsets of $\Lambda$ and 
$\mu_\La$ is the  grand canonical Gibbs measure on $\Gamma_\Lambda$ corresponding to $\phi$ and empty boundary condition. 
 The problems are: 1) prove that the generator of $\EE_\La$ has a spectral gap, 2) estimate how quickly it  shrinks as $\La\nearrow\X$.
 Note that problems of such type have been  studied in the lattice case, see e.g.\ \cite{CM}.

Finally, in our forthcoming paper \cite{KKL}, we discuss another possible scaling of the family of the Kawasaki dynamics \eqref{jhvhgjvf} with $\varepsilon =1$ which leads to a family of infinite particle diffusions. The latter family includes as a special case the gradient dynamics, see e.g.\ \cite{Fritz,AKR2} and  the refeferences therein.

  \begin{center}
{\bf Acknowledgements}\end{center}

 The authors acknowledge the financial support of the SFB 701 `` Spectral
structures and topological methods in mathematics'', Bielefeld University. We would like to thank Zeev Sobol for useful discussions.

\end{document}